\documentclass{amsart}
\usepackage{amsmath,amssymb,amsthm,amscd, mathtools, latexsym}
\usepackage{enumerate,varioref, dsfont, fancyhdr}
\usepackage{tikz-cd}
\usepackage{enumitem}
\usepackage{multicol}
\usepackage{hyperref}
\usepackage{url}
\usepackage{textcomp}

\newtheorem{Thm}{Theorem}[section]
\newtheorem{Lem}{Lemma}[section]
\newtheorem{Prop}{Proposition}[section]

\newtheorem{Cor}{Corollary}[section]

\newtheorem{Rem}{Remark}[section]
\newtheorem{Obs}{Observation}[section]
\newtheorem{Def}{Definition}[section]
\newtheorem{claim}{Claim}[section]

{\theoremstyle{plain}

}

\def\rit{{\mathbb R}}
\def\cit{{\mathbb C}}
\def\nit{{\mathbb N}}
\def\qit{{\mathbb Q}}
\def\zit{{\mathbb Z}}

\def\pit{{\mathbb P}}

\def\0{{\mathcal O}}

\def\A{{\mathfrak A}}

\def\I{{\mathcal I}}

\def\C{{\mathcal C}}

\def\M{{\mathcal M}}

\def\A{{\mathcal A}}
\def\E{{\mathcal E}}
\def\W{{\mathcal W}}

\title{The Chow group mod $\ell$ of a product of elliptic curves}
\author{Humberto A. Diaz}
\date{}
\begin{document}
\begin{abstract}
Generalizing work of Schoen, we prove that the Chow group modulo $\ell$ of a product of $3$ or more very general complex elliptic curves is infinite. 
\end{abstract}
\maketitle

\section{Introduction}
A recent result of Totaro \cite{T} shows that for any prime $\ell$ the Chow group mod $\ell$ of a very general Abelian threefold is infinite, generalizing a result of Rosenschon and Srinivas in \cite{RS}. Our goal will be to prove the following analgous result, which extends work of Schoen in \cite{S1}:
\begin{Thm}\label{1} Let $E_{1}$, $E_{2}, \ldots E_{n}$ be very general complex elliptic curves. Then,
$$CH^{d} (E_{1}\times E_{2} \times \dots \times E_{n}) \otimes \zit/\ell$$
is infinite for all primes $\ell$ and $2 \leq d \leq n-1$.
\end{Thm}
\noindent The key case is when $d=2$ and $n=3$ since the other cases may be obtained by pulling back along the projection to the triple product and then intersecting with a symmetric divisor. 
Our strategy in this case will be as follows. As noted by Totaro in \cite{T}, it suffices to prove the corresponding result modulo $\ell^{r}$ for $r >>0$. To this latter end, we first find a cycle which does not vanish modulo $\ell^{r}$. This will be an ersatz Ceresa cycle $\gamma$ which is a spreading out of the cycle considered in \cite{Bl} and \cite{BST} to a certain $5$-dimensional family of products of elliptic curves. It turns out that (using some rather general methods of Voisin in \cite{V}) one can show that the geometry of the total space of this family is rather simple; this allows for a self-contained geometric proof of the non-vanishing of the normal function of $\gamma$ (in singular cohomology). More precisely, $\sim_{rat} = \sim_{hom}$ for this family so that proving homological non-triviality for $\gamma$ reduces to proving rational non-triviality for $\gamma$. It seems likely that this approach might work for other families of varieties for which the methods of \cite{V} apply. This approach, however, does not work for the usual Ceresa cycle of genus $3$ since this cycle is not defined on the universal Jacobian of genus $3$ (with level $N$-structure), $J(N) \to \M_{3} (N)$ but only over some cover of the generic fiber. (In this case, the authors of \cite{RS} and \cite{T} invoke a rather sophisticated result of Hain in \cite{H} to prove the non-triviality of the normal function.) The results of \cite{BE} then give the required non-triviality mod $\ell^{r}$. In order to obtain infinitely many such cycles, we will modify Nori's usual isogeny argument (using his suggestion in \cite{N}).  

\subsection*{Acknowledgements}

The author would like to thank Burt Totaro for his interest in early drafts of this paper. He would also like to thank Chad Schoen for some useful comments and Robert Laterveer for calling his attention to the above-mentioned method of Claire Voisin.

\section{Construction of the main cycle}
Let $k\subset \cit$ be algebraically closed and $U \subset \pit(H^{0} (\pit^{2}, \mathcal{O}_{\pit^{2}}(2))$ be an open subset consisting only of smooth conics. Then, we consider the universal family of conics over $U$:
\[ \begin{CD}  \mathcal{Q} @>>> & \pit^{2} \\ @VhVV & \\ U &  \end{CD}\]
There is a finite morphism $\phi: \pit^{2} \to \pit^{2}$ defined by
\[ [x, y, z] \mapsto [x^{2}, y^{2}, z^{2}] \]
whose Galois group is $G:= \zit/2 \times \zit/2$ and is generated by the involutions:
\[ [x, y, z] \xmapsto{\sigma_{1}} [-x, y, z], \ [x, y, z] \xmapsto{\sigma_{2}} [x, -y, z]  \]
and denote by $\sigma_{3} := \sigma_{1}\sigma_{2}$. We then consider the Cartesian product:
\begin{equation}\label{cart-diag} \begin{CD}  \mathcal{C} @>q>> &\mathcal{Q}  \\ @VVV & @VVV\\ \pit^{2} @>\phi>> & \pit^{2} \end{CD}\end{equation}
There is a natural morphism $f: \mathcal{C} \xrightarrow{} U$, and after possibly shrinking $U$, we may assume that all the fibers of $f$ are smooth. The generic fiber is then a genus $3$ curve in $\pit^{2}_{k(U)}$ defined by the following equation:
\[ ax_{1}^{4}+bx_{2}^{4}+cx_{3}^{4}+dx_{1}^{2}x_{2}^{2}+ex_{2}^{2}x_{3}^{2}+fx_{2}^{2}x_{3}^{2}=0  \]
The involutions above then induce an action on $\mathcal{C} \xrightarrow{f} U$, which gives rise to the quotients:
\[ p_{j}: \C \to \E_{j}:= \mathcal{C}/\sigma_{j}, \ g_{j}: \mathcal{E}_{j} \to U, \ j=1,2,3 \]
where $g_{j} \circ p_{j} = f$. Note then that on $\pit^{2}$ the involution $\sigma_{j}$ fixes the line defined by $x_{j} = 0$ (as well as an isolated point). Again, upon shrinking $U$, we may assume that none of the fibers of $f$ pass through any of these isolated points. Thus, $\sigma_{j}$ fixes precisely 4 points on the general fiber of $f$. A Riemann-Hurwitz argument then shows that the fibers of $g_{j}$ all have genus 1. Moreover, there are induced quotient maps:
\[ q_{j} : \E_{j} \to \mathcal{Q} = \C/G \]
as well as involutions, $(-1)_{\E_{j}}: \E_{j} \to \E_{j},$ which induce the usual action by $-1$ on the geometric generic fiber of $g_{j}$. Then, we consider the fiber product
\[ g_{1} \times_{U} g_{2}  \times_{U} g_{3}: \E_{1} \times_{U} \E_{2} \times_{U} \E_{3} \to U \]
and use the short-hand $g:= g_{1} \times_{U} g_{2}  \times_{U} g_{3}$ and $\A:= \E_{1} \times_{U} \E_{2} \times_{U} \E_{3}$. There is then an action by \[\tilde{G} = \zit/2 \times \zit/2 \times \zit/2\] on $g: \A \to U$ generated by involutions
\[ \tilde{\sigma}_{j} : \A \to \A, \ j=1,2,3 \]
which act by $(-1)_{\E_{j}}$ on the $j^{th}$ factor and the identity on the remaining factors. We also use the notation:
\[ \tilde{\sigma}_{jk} := \tilde{\sigma}_{j}\circ\tilde{\sigma}_{k}, \ \tilde{\sigma}_{123}:= \tilde{\sigma}_{1}\circ\tilde{\sigma}_{2}\circ\tilde{\sigma}_{3} \]
Moreover, there is a morphism:
\[ \rho: \C \xrightarrow{\Delta^{(3)}} \C^{(3)} := \C \times_{U} \C \times_{U} \C \xrightarrow{p} \A \]
where $p:= p_{1} \times_{U} p_{2} \times_{U} p_{3}$ and 
\[ \rho_{jk}: \C \xrightarrow{\rho} \A \xrightarrow{\pi_{jk}} \E_{j} \times_{U} \E_{k} \]
 Then, we have the following lemma:
\begin{Lem}\label{basics}
The generic fiber of $g: \A \to U$ is geometrically isogenous to the Picard variety of the generic fiber of $f: \C \to U$ and there is a commutative diagram:
\[\begin{CD} \C @>{\sigma_{i}}>> & \C \\ @V{\rho}VV & @V{\rho}VV \\ \A @>{\tilde{\sigma}_{jk}}>> & \A \end{CD} \]
for $i \not\in \{ j, k \}$. Further, $\rho$ is an imbedding; moreover, $\rho_{jk}$ is an immersion and is injective away from the fixed points of $\sigma_{i}$ for $i \not\in \{ j, k \}$.
\begin{proof} The second statement follows directly from the constructions of $\sigma_{i}$ and $\tilde{\sigma}_{jk}$. For the third statement, suppose $\rho(x) = \rho(x')$ for $x \neq x'$; i.e., that
\begin{equation}\sigma_{1}(x) = x', \ \sigma_{2}(x) = x', \ \sigma_{3}(x) = x' \label{3}\end{equation}
But then one has $x' = \sigma_{3}(x) = \sigma_{1}(\sigma_{2}(x)) = \sigma_{1}(x') = x$, which is a contradiction. For the fourth statement, the immersion part follows from noting that
\[ d\rho_{jk} = (dp_{j}, dp_{k}) : TC_{x} \to TE_{j, x} \oplus TE_{k, x} \]
and this vanishes only where $dp_{j}$ and $dp_{k}$ vanish simultaneously, which does not happen. Indeed, possibly after shrinking $U$, we can assume that none of the fibers of $\C \to U$ contain fixed points of both $\sigma_{j}$ and $\sigma_{k}$. For generic injectivity, note that if only (\ref{3}) holds for only 2 out of the 3 involutions, $x$ must be a fixed point of the remaining involution. Finally for the first statement, let $x \in U$, $C_{x}$ and $A_{x}$ be fibers of $f$ and $g$, respectively. Then, using the universal property of the Albanese, there is a commutative diagram:
\[ \begin{tikzcd} C_{x} \arrow[hook]{r} \arrow{rd}{\rho} & Pic^{0}(C_{x})\arrow{d}\\ & A_{x} \end{tikzcd} \] 
So, it suffices to show that the map $\rho_{x}^{*}: H^{0} (A_{x}, \Omega^{1}_{A_x}) \to H^{0} (C_{x}, \Omega^{1}_{C_x})$ is an isomorphism. However, since both sides have the same dimension, we are reduced to proving injectivity. To this end, $A_{x} = E_{1,x} \times E_{2,x} \times E_{3,x}$, so that \[ H^{0} (A_{x}, \Omega^{1}_{A_x}) = \pi_{1}^{*}H^{0} (E_{1,x}, \Omega^{1}_{E_{1,x}}) \oplus \pi_{2}^{*}H^{0} (E_{2,x}, \Omega^{1}_{E_{2,x}}) \oplus \pi_{3}^{*}H^{0} (E_{3,x}, \Omega^{1}_{E_{3,x}}) \]
Also, $\pi_{j}\circ\rho_{x} = p_{j}: C_{x} \to E_{j,x}$. So, injectivity will follow from the statement that
\begin{equation} p_{j}^{*}H^{0} (E_{j,x}, \Omega^{1}_{E_{j,x}}) \cap p_{k}^{*}H^{0} (E_{k,x}, \Omega^{1}_{E_{k,x}}) =0\label{inter}\end{equation}
for $j \neq k$. For this, observe that the left hand side of (\ref{inter}) is invariant under $\sigma_{j}$ and $\sigma_{k}$ and, hence, under $G$. However, $H^{0} (C_{x}, \Omega^{1}_{C_x})^{G} = H^{0} (Q_{x}, \Omega^{1}_{Q_x}) = 0$, which gives the desired equality.
\end{proof}
\end{Lem}
The main cycle is then defined as an analogue of the Ceresa cycle:
\begin{equation} \boxed{\gamma := \rho_{*}([\C]) - \tilde{\sigma}_{123}^{*}\rho_{*}([\C]) \in CH^{2} (\A)} \end{equation}
We observe that $\gamma$ is anti-invariant under $\tilde{\sigma}_{123}$ and invariant under $\tilde{\sigma}_{jk}$ by Lemma \ref{basics} and, hence, is anti-invariant with respect to each $\tilde{\sigma}_{j}$. 
\section{The normal function of $\gamma$}\label{normal}

We now take $k = \cit$ and let $g_{j}: \E_{j} \to U$ and $g: \A \to U$ be as in the previous section. Then, there is a decomposition of local systems over $U$ given by:
\[ R^{1}g_{*}\qit \cong \bigoplus_{j=1}^{3} R^{1}g_{j*}\qit \]
Since $\bigwedge^{2} R^{1}g_{j*}\qit \cong \qit$, we also have
\[ R^{3}g_{*}\qit \cong \bigwedge^{3}R^{1}g_{*}\qit \cong \bigoplus_{j=1}^{3} R^{1}g_{j*}\qit \oplus PR^{3}g_{*}\qit \]
where $PR^{3}g_{*}\qit = R^{1}g_{1*}\qit\otimes R^{1}g_{2*}\qit\otimes R^{1}g_{3*}\qit$. (The notation is suggestive of primitive cohomology, which is not an irreducible local system in the case at hand; so, $PR^{3}g_{*}\qit$ will be considered in its place.) Now, we can find a relative correspondence:
\[ P \in CH^{3} (\A \times_{U} \A) \]
for which $P_{*}(Rg_{*}\qit) \cong PR^{3}g_{*}\qit [-3] \in D^{b} (S, \qit)$, the derived category of $\qit$-local systems over $U$. Indeed, we can select
\[ \pi_{1,j} \in CH^{1} (\E_{j} \times_{U} \E_{j}) \]
for which $\pi_{1,j*}(Rg_{j*}\qit) = R^{1}g_{j*}\qit[-1]$. (One can set $\pi_{1,j} = \Delta_{\E_{j}} - \Gamma_{\tilde{\sigma}_{j}}$ and this has the desired effect.) Then set
\[ P = \pi_{1,1} \times_{U} \pi_{1,2} \times_{U} \pi_{1,3} \in CH^{3} (\A \times_{U} \A) \]
 With this definition, we have:
\begin{Obs}\label{important} $P_{*}CH^{j} (\A)_{\qit}$ is the subspace of $CH^{j} (\A)_{\qit}$ on which $\tilde{\sigma}_{j}$ acts by $-id$.
\end{Obs}
\noindent Further, since $P_{*}(Rg_{*}\qit) \cong PR^{3}g_{*}\qit [-3]$, we can apply $H^{4}$ to obtain:
\[P_{*}H^{4} (\A, \qit) \cong H^{1} (U, PR^{3}g_{*}\qit) \]
\begin{Prop}\label{key} $P_{*}\gamma \neq 0 \in H^{1} (U, PR^{3}g_{*}\qit)$

\begin{proof} We will need to show that $P_{*}\gamma \neq 0 \in H^{4} (\A, \qit)$; an important step towards this is the following lemma, which shows that the geometry of $\A$ is rather simple.
\begin{Lem} \begin{enumerate}[label=(\alph*)]
\item\label{geo} $\C^{(3)} := \C \times_{U} \C \times_{U} \C$ contains $C^{(3)}_{0}$, an open subset of a projective bundle over the blow-up of $\pit^{2} \times \pit^{2}\times \pit^{2}$ along a union of (transversely-intersecting) rational varieties.
\item\label{cyc} The cycle class map $CH^{2} (\A)_{\qit} \to H^{4} (\A, \qit)$ is injective.
\end{enumerate}
\begin{proof} For \ref{geo}, we adapt an argument of Voisin (Proposition 3.11 in \cite{V}), although the case at hand is less technical. We note that the universal conic over $U \subset H^{0} (\pit^{2}, \mathcal{O}_{\pit^{2}}(2))$ is given by
\[ \mathcal{Q} = \{ (Q, x) \in U \times \pit^{2} \ | \ Q \in U, \ Q(x) = 0  \} \]
We set $(\pit^{2})^{3} := \pit^{2} \times \pit^{2} \times \pit^{2}$ and consider the triple fiber product over $U$:
\[ \mathcal{Q}^{(3)} := \mathcal{Q} \times_{U} \mathcal{Q} \times_{U} \mathcal{Q} = \{ (Q, x, y, z) \in U \times (\pit^{2})^{3} \ | \ Q \in U, \ Q(x) = Q(y) = Q(z) = 0  \} \]
There is the natural map $p^{(3)}: \mathcal{Q}^{(3)} \to (\pit^{2})^{3}$ defined as the projection onto the last 3 factors. Now, we consider the closed subvarieties $\Delta_{ij} \subset \pit^{(3)}$ defined as the set of $(x_{1}, x_{2}, x_{3}) \in (\pit^{2})^{3}$ such that $x_{i} = x_{j}$ and the blow-up $\epsilon: \widetilde{(\pit^{2})^{3}} \to (\pit^{2})^{3}$ along
\[ \Delta_{12} \cup \Delta_{23} \cup \Delta_{13}\] 
Since $\mathcal{Q} \to U$ is a relative curve, we see that $(p^{(3)})^{-1} (\Delta_{ij})$ is of codimension $1$ and thus by the universal property of blow-up $p^{(3)}: \mathcal{Q}^{(3)} \to (\pit^{2})^{3}$ factors through $\epsilon: \widetilde{(\pit^{2})^{3}} \to \pit^{(3)}$, giving a morphism $\tilde{p}^{(3)}: \mathcal{Q}^{(3)} \to \widetilde{(\pit^{2})^{3}}$. Now, we let 
\[ {\mathcal{Q}}^{(3)}_{0} = {\mathcal{Q}}^{(3)} \setminus \Delta_{\mathcal{Q}}^{(3)}, \ \widetilde{(\pit^{2})^{3}}_{0} = \widetilde{(\pit^{2})^{3}} \setminus \epsilon^{-1}(\Delta_{\pit^{2}}^{(3)}) \]
where $\Delta_{\mathcal{Q}}^{(3)}$ and $\Delta_{\pit^{2}}^{(3)}$ denote the small diagonals in $\mathcal{Q}^{(3)}$ and $(\pit^{2})^{3}$, respectively. By construction, $\widetilde{(\pit^{2})^{3}}_{0}$ consists of length $3$ subschemes of $\pit^{2}$ which are supported on at least $2$ distinct points. Then, there is a map
\[ {p}^{(3)}_{0}: {\mathcal{Q}}^{(3)}_{0} \to \widetilde{(\pit^{2})^{3}}_{0}\]
such that for any length $3$ subscheme $Z \in \widetilde{(\pit^{2})^{3}}_{0}$ the fiber $({p}^{(3)}_{0})^{-1} (Z)$ consists of all conics $Q \in U$ such that $Q$ vanishes on $Z$. It is an easy result that such a subscheme imposes $3$ independent conditions on conics, so that the dimension of any fiber $({p}^{(3)}_{0})^{-1} (Z)$ is an open subset of $\pit^{2}$. In fact, ${\mathcal{Q}}^{(3)}_{0}$ is an open subscheme of a projective bundle over $\widetilde{(\pit^{2})^{3}}_{0}$, for which the fiber over $Z \in \widetilde{(\pit^{2})^{3}}_{0}$ is $H^{0} (\pit^{2}, \I_{Z} (2))$.\\
\indent Recall the morphism $\phi: \pit^{2} \to \pit^{2}$ defined in the previous section by $[x,y,z] \mapsto [x^{2}, y^{2}, z^{2}]$ and now consider the triple product of this map, $\phi^{3}: (\pit^{2})^{3} \to (\pit^{2})^{3}$. Now, $(\phi^{3})^{-1} (\Delta_{jk})$ consists of $(x_{1}, x_{2}, x_{3}) \in (\pit^{2})^{3}$ satisfying \[x_{j}^{2}-x_{k}^{2} = (x_{j} - x_{k})(x_{j} + x_{k}) = 0\] So, we let $\hat{\epsilon}: \widehat{(\pit^{2})^{3}} \to (\pit^{2})^{3}$ be the blow-up of $(\pit^{2})^{3}$ along the subvariety 
\[ (\phi^{3})^{-1} (\Delta_{12}) \cup (\phi^{3})^{-1} (\Delta_{23}) \cup (\phi^{3})^{-1} (\Delta_{13}) \]
We then obtain a morphism $\tilde{\phi}^{3}: \widehat{(\pit^{2})^{3}} \to \widetilde{(\pit^{2})^{3}}$. From the Cartesian diagram (\ref{cart-diag}), it is straightforward to see that there is a Cartesian diagram:
\[\begin{CD} \C^{(3)} @>{q^{(3)}}>> \mathcal{Q}^{(3)}\\
@VVV @V{p^{(3)}}VV\\  
\widehat{(\pit^{2})^{3}} @>{\tilde{\phi}^{3}}>> \widetilde{(\pit^{2})^{3}}\end{CD}\]
Thus, the open set $\C^{(3)}_{0}:=(q^{(3)})^{-1}(\mathcal{Q}^{(3)}_{0}) \subset \C^{(3)}$ is an open subset of a projective bundle over $\widehat{(\pit^{2})^{3}}_{0}:= (\tilde{\phi}^{3})^{-1}(\widetilde{(\pit^{2})^{3}}_{0}) \subset \widehat{(\pit^{2})^{3}}$, which gives the proof of \ref{geo}.\\
\indent For \ref{cyc}, by applying Lemmas 2.1, 2.2, and 2.3 of \cite{V} to $\C^{(3)}_{0}$, it follows that the cycle class map:
\[CH^{2}(\C^{(3)}_{0})_{\qit} \to H^{4}(\C^{(3)}_{0}, \qit) \]
injective. Since the complement of $\C^{(3)}_{0}$ in $\C^{(3)}$ has codimension $\geq 2$, a Gysin exact sequence argument shows that 
\begin{equation} CH^{2}(\C^{(3)})_{\qit} \to H^{4}(\C^{(3)}, \qit)\label{pre-inv} \end{equation}
is also injective. Finally, since $\A = \E_{1} \times_{U} \E_{2} \times_{U} \E_{3}$ is a quotient of $\C^{(3)}$ by $H = \zit/2 \times \zit/2 \times \zit/2$, applying the functor $(\phantom{-})^{H}$ to (\ref{pre-inv}), we obtain the statement of the lemma. 
 
\end{proof}
\end{Lem}
Now, suppose by way of contradiction that $P_{*}\gamma = 0 \in H^{4} (\A, \qit)$; then, by the previous lemma, we must have $P_{*}\gamma = 0 \in CH^{2} (\A)_{\qit}$. Since $\gamma$ is anti-invariant under $\tilde{\sigma}_{j}$, it follows by Observation \ref{important} that $\gamma = 0 \in CH^{2} (\A)_{\qit}$ and, hence, by specialization:
\begin{equation} \gamma_{x} = 0 \in CH^{2} (E_{1,x} \times E_{2,x} \times E_{3,x})_{\qit}\label{vanishing} \end{equation}
for each $x \in U$. Now, since $x$ varies in a $5$-dimensional family, we may fix $E_{2, x} = E$ and $E_{3, x}= E'$ and allow $E_{1, x} = E_{t}$ to vary in a $1$-parameter family indexed by $t$.
The desired contradiction is obtained using the following lemma:
\begin{Lem} If (\ref{vanishing}) holds, then $CH^{2} (K)_{\qit}$ has rank $1$, where $K$ is the Kummer surface of $E \times E'$.
\begin{proof}[Proof of Lemma] By Lemma \ref{basics}, one has:
\[C_{t} \subset E_{t} \times E \times E' \]
and the projection to the last two factors gives an immersed curve $\rho_{23}(C_{t}) \subset E \times E'$ whose only singularities are nodes contained in the set of $2$-torsion points of $E \times E'$ (i.e., the fixed points of $\tilde{\sigma}_{23}$). Now, we let $\widetilde{E \times E'}$ denote the blow-up of $E \times E'$ along the $2$-torsion points. There is a lift $C_{t} \subset E_{t} \times \widetilde{E \times E'}$ for which $\rho_{23}(C_{t})$ is smooth (since its only singular points are nodes in the set of blown-up points). The action of $\tilde{\sigma}_{23} = (-1)_{E \times E'}$ lifts to $\widetilde{E \times E'}$ and the corresponding action on $E_{t} \times \widetilde{E \times E'}$ induces the action of $\sigma_{1}$ on $C_{t}$ by Lemma \ref{basics}. Taking the quotient by $\sigma_{1}$ gives
\[ E_{t} \subset E_{t} \times K \]
where $K$ is the Kummer surface of $\widetilde{E \times E'}$. By construction, the projection of $E_{t}$ onto $K$ is an isomorphism, and so we have an inclusion $j_{t} : E_{t} \hookrightarrow K$. Now, since $K$ is a $K3$ surface, $|j_{t}(E_{t})|$ is a $1$-dimensional linear system and induces the structure of an elliptic surface on $K$. In fact, we could have chosen the $1$-parameter family indexed by $t$ so that $|j_{t}(E_{t})|$ does not depend on $t$. Thus, $|j_{t}(E_{t})|$ induces the structure of an elliptic surface on $K$. Finally, the action of $\tilde{\sigma}_{1}$ on $E_{t} \times E \times E'$ descends to $E_{t} \times K$, acting as $(-1)_{E_{t}}$ on $E_{t}$ and the identity on $K$.\\
\indent By our assumption, the main cycle $\gamma_{t} \in CH^{2} (E_{t} \times E \times E')_{\qit}$ vanishes. We note that
\[ \gamma_{t} = \rho_{*}(C_{t}) - \tilde{\sigma}_{123}^{*}\rho_{*}(C_{t}) = \rho_{*}(C_{t}) - \tilde{\sigma}_{1}^{*}\rho_{*}(C_{t}) \]
since $\rho(C_{t})$ is $\tilde{\sigma}_{23}$-invariant. Pulling $\gamma_{t}$ back via the blow-up and pushing it forward to $E_{t} \times K$, we obtain 
\[ \overline{\gamma}_{t} = (id \times j_{t}){*}(\Delta_{E_{t}}) - ((-1)_{E_{t}} \times j_{t})_{*}(\Delta_{E_{t}}) = (id \times j_{t})_{*} (\Delta_{E_{t}} - \Gamma_{(-1)_{E_{t}}}) \in CH^{2} (E_{t} \times K)_{\qit}\]
Now, we may view $\overline{\gamma}_{t}$ as a correspondence between $E_{t}$ and $K$. By Liebermann's lemma, we have
\[ \overline{\gamma}_{t} = (id \times j_{t})_{*} (\Delta_{E_{t}} - \Gamma_{(-1)_{E_{t}}}) = \Gamma_{j_{t}}\circ(\Delta_{E_{t}} - \Gamma_{(-1)_{E_{t}}}) \in Cor^{1} (E_{t} \times K)\]
Let $e_{t} \in CH^{1} (E_{t})_{\qit}$ denote the class of the identity element of $E_{t}$; then, we have the decomposition:
\[ CH^{1}(E_{t})_{\qit} = \qit\cdot e_{t} \oplus CH^{1}_{(1)} (E_{t})_{\qit} \]
where $CH^{1}_{(1)} (E_{t})_{\qit}$ is the summand on which $(-1)_{E_{t}}$ acts by $-id$. Since $\gamma_{t} = 0$ for all $t$, it follows that for any $\alpha \in CH^{1} (E_{t})_{\qit}$ of degree $d$ we have
\[ 0 = j_{t*}(\alpha) =j_{t*}(\alpha - de_{t})  \]
from which it follows that $j_{t*}(CH^{1}(E_{t})_{\qit}) = \qit j_{t*}(e_{t}) \subset CH^{2} (K)$. Moreover, $j_{t*}(e_{t}) = j_{s*}(e_{s})$ for all $s, t$ (since both lie on a rational curve). Thus, the image of the direct sum of pushforwards:
\[ \oplus_{t} CH^{1} (E_{t})_{\qit} \xrightarrow{\oplus j_{t*}} CH^{2} (K)_{\qit} \]
has rank $1$. By our construction, $K$ is an elliptic surface whose fibers are $E_{t}$ and the $E_{t}$ account for all but finitely many of the fibers. This would imply that $CH^{2} (K)_{\qit}$ has rank $1$, as desired.
\end{proof}
\end{Lem}
\end{proof}
\end{Prop}
\begin{Cor}\label{main} Let $k_{0} \subset \cit$ be an algebraically closed field such that $k_{0} (U) \subset \cit$ and let $G_{k_{0}(U)}$ be the absolute Galois group and 
\[ P_{*}\gamma \in CH^{2} (A_{\overline{k_{0}(U)}}) \]
be the cycle constructed above (on the geometric generic fiber of $\A \to U$). Then, for each prime $\ell$ there is some $r >>0$ for which the cycle class
\begin{equation} P_{*}\gamma \neq 0 \in H^{1} (G_{k_{0}(U)}, PH^{3}_{\text{\'et}} (A_{\overline{k_{0}(U)}}, \zit/\ell^{r}(2)))\label{Pgamma} \end{equation}
Moreover, the class (\ref{Pgamma}) is invariant with respect to $\tilde{\sigma}_{23}$ and anti-invariant with respect to $\tilde{\sigma}_{123}$.
\begin{proof} From the Gysin sequence, the natural map
\[H^{1}_{\text{\'et}} (U_{k_{0}}, PH^{3}_{\text{\'et}} (A_{\overline{k_{0}(U)}}, \zit/\ell^{r}(2))) \to H^{1} (G_{k_{0}(U)}, PH^{3}_{\text{\'et}} (A_{\overline{k_{0}(U)}}, \zit/\ell^{r}(2))) \]
is injective. Moreover, since \'etale cohomology is unchanged under taking algebraically closed extensions, it suffices to show that 
\[ P_{*}\gamma \neq 0 \in H^{1}_{\text{\'et}} (U_{\cit}, PH^{3}_{\text{\'et}} (A_{\overline{\cit(U)}}, \zit/\ell^{r}(2)))\] 
Then, using the compatibility isomorphism with singular cohomology, it suffices to show that this is true for singular cohomology. On the other hand, from the above proposition, we have that
\[ P_{*}\gamma \neq 0 \in H^{1} (U_{\cit}, PH^{3} (A_{{\cit}}, \qit(2))) \] 
which implies $ P_{*}\gamma \in H^{1} (U_{\cit}, PH^{3} (A_{{\cit}}, \zit(2)))$ is non-torsion, which means that \[ P_{*}\gamma \neq 0 \in H^{1} (U_{\cit}, PH^{3} (A_{{\cit}}, \zit(2))) \otimes \zit/\ell^{r}\]
for $r >>0$. Since $H^{3} (A_{{\cit}}, \zit)$ is torsion-free, the universal coefficient theorem then gives the desired result.

\end{proof}
\end{Cor} 
\section{Nontriviality of $\gamma$}

We would like to realize the above construction as a cycle on the general product of three elliptic curves. For this, we will need the following result on the period map associated to $\A \to U$. 

\begin{Lem} Let $k \subset \cit$ be an algebraically closed field $X(1)$ the $j$-line over $k$. Then, the period map $U \to X(1)^{\times 3}$ induced by $g: \A \to U$ is dominant.
\begin{proof} Indeed, note that by Lemma \ref{basics} $A_{x}$ is isogenous to $Pic^{0}(C_{x})$, so it suffices to prove the corresponding statement for the period map $U \to \A_{3}$ induced by $Pic^{0} (\C) \to U$ (where $A_{3}$ is the coarse moduli space of ppav's of dimension $3$). But by the Torelli theorem, it suffices to prove the statement for the map $U \to \M_{3}$ (induced by $C \to U$), which will be done once we show that the image of $U \to \M_{3}$ has dimension $3$ or, equivalently, that the fibers are generally of dimension $2$. This latter is precisely the claim below:
\begin{claim} For the general $x \in U$, $\{ y \in U \  | \ C_{y} \cong C_{x} \}$
has dimension $2$.
\begin{proof}[Proof of Claim] First note that an isomorphism $C_{y} \cong C_{x}$ is induced by an automorphism $\phi: \pit^{2} \to \pit^{2}$ that maps $C_{y}$ to $C_{x}$. We would like to show that $\phi$ must lie in an algebraic group of dimension $2$. To this end, note that for the general $x \in U$, $Aut(C_{x}) = G$, (one can see this by specializing to the Fermat quartic curve whose automorphism group is the wreath product $\zit/4 \wr \S_{3}$ and then noting that only the elements of $G$ spread out to automorphisms of the general fiber of $\C \to U$). Thus, such a $\phi \in PGL_{3}(k)$ must lie in $N(G)$, the normalizer of $G \leq PGL_{3}(k)$. Then, using a dimension argument, it is easy to see that $k[G]^{\times} \hookrightarrow GL_{3}(k)$ is the group of diagonal matrices $D$, so that $N(G) = N(D)$. However, since $D$ is a torus over $k$, it is a well-known fact that the identity component of $N(D)$ is the centralizer of $D$, $C(D) = D$. Thus, $N(D)$ has dimension $2$, as desired.
\end{proof}
\end{claim}
\end{proof}
\end{Lem}
\begin{Cor}\label{rigid} Let $k_{0} \subset \cit$ be an algebraically closed field for which $k_{0} (U) \subset \cit$. Then, $A_{k_{0}(U)} \times_{k_{0}(U)} \cit$ is isomorphic as a scheme to the very general product of $3$ complex elliptic curves, $E_{1} \times E_{2} \times E_{3}$. Moreover, for each prime $\ell$ and $r \in \nit$, we have
\[ CH^{2} (A_{\overline{k_{0}(U)}}) \otimes \zit/\ell^{r} \cong CH^{2} (E_{1} \times E_{2} \times E_{3}) \otimes \zit/\ell^{r} \]
\begin{proof} The first statement follows from the previous lemma and the second statement follows from Suslin rigidity \cite{Su} in the form of \cite{Le}.
\end{proof}
\end{Cor}
In order to prove the nontriviality of $\gamma$, we will need the following lemma, whose proof relies on a deep result of Bloch-Esnault in \cite{BE} as well as the Merkurjev-Suslin theorem. Similar results for more general Abelian varieties have been obtained by Schoen in \cite{S2}, \cite{S1} and Totaro in \cite{T}. 
\begin{Lem} For each prime $\ell$ and $r \in \nit$ and very general complex elliptic curves $E_{1}$, $E_{2}$ and $E_{3}$, 
\[ P_{*}CH^{2} (E_{1} \times E_{2} \times E_{3})[\ell^{r}] = 0 \]
\begin{proof} While general results of this type are well-known to the experts, it is possible to determine explicitly the $\ell^{r}$-torsion in this case. We refer the reader to the appendix for this.
\end{proof}
\end{Lem}
\begin{Cor}\label{main-2} Let $\ell$ be a prime, $m \in \nit$ and $E_{1}$, $E_{2}$, $E_{3}$ be very general elliptic curves over $\cit$. Then, there exists $\gamma \in CH^{2} (E_{1} \times E_{2} \times E_{3})$ for which  \[ \gamma \neq 0 \in CH^{2} (E_{1} \times E_{2} \times E_{3}) \otimes \zit/\ell^{r}\] is anti-invariant with respect to $(-1)_{E_{j}}$ for $r >>0$.
\begin{proof} By Corollary \ref{rigid}, we have
\[ CH^{2} (A_{\overline{k_{0}(U)}}) \otimes \zit/\ell^{r} \cong CH^{2} (E_{1} \times E_{2} \times E_{3}) \otimes \zit/\ell^{r} \]
so it suffices to prove the corresponding result for $A_{\overline{k_{0}(U)}}$. Suppose by way of contradiction that for $r>>0$, there is some $\delta \in CH^{2} (A_{\overline{k_{0}(U)}})$ such that
\[ \gamma = \ell^{r}\cdot\delta \]
Then, following the argument of \cite{BE} Proposition 4.1 (or \cite{S2} Proposition 6.2 or \cite{T} \S 2), we would like to show that $P_{*}\delta$ is defined over $k_{0}(U)$. To this end, let $G_{k_{0}(U)}$ be the absolute Galois group and $g \in G_{k_{0}(U)}$. Since $\gamma$ is defined over $k_{0}(U)$, we have $g(\gamma) = \gamma$ for all $g \in G_{k_{0}(U)}$. Thus,
\[ 0 = m(\gamma - g(\gamma)) = \ell^{r}(\delta - g(\delta)) \]
So, $\delta - g(\delta)$ is $\ell^{r}$-torsion, which means that by the previous lemma, 
\[ 0=P_{*}(\delta - g(\delta)) = P_{*}\delta - g(P_{*}\delta)) \]
So, $P_{*}\delta$ descends to a cycle over $k_{0}(U)$ which implies that
\[ P_{*}\gamma = 0 \in CH^{2} (A_{{k_{0}(U)}}) \]
Applying the cycle class map, we obtain a contradiction to Corollary \ref{main}. As in \cite{T}, the exact same argument applied to $2\cdot\gamma$ shows that $\gamma$ satifies the required anti-invariance property.
\end{proof}
\end{Cor}

\section{Proof of Theorem \ref{1}}

Beginning with $\gamma$, we will use Nori's isogeny argument in \cite{N} to obtain the required infinitely many cycles mod $\ell^{r}$. We let $N \geq 3$ be a positive integer such that $\ell \nmid N$ and let $X(N)$ be the fine moduli space of elliptic curves with full level-$N$ structure, let $\E(N) \to X(N)$ be the universal elliptic curve. Then, consider the triple product
\begin{equation} \A(N) := \E(N) \times \E(N) \times \E(N) \to W(N):= X(N) \times X(N) \times X(N) \label{universal} \end{equation}
Observe that $W(N)$ is the fine moduli space parametrizing products of $3$ elliptic curves with full level $N$-structure and $\A(N)$ is the universal triple roduct.
Further, denote by 
\[ K_{N} = \cit(W(N)) \]
the function field of $W(N)$ and let $A_{K_{N}}$ be the generic fiber of (\ref{universal}). Note that there are $3$ involutions, $\tilde{\sigma}_{1}$, $\tilde{\sigma}_{2}$, and $\tilde{\sigma}_{3}$, on $A_{K_{N}}$ (induced by the action of $-1$ on each of the three factors). Then, by Corollary \ref{main-2}, there exists $\gamma \in CH^{2} (A_{\overline{K_{N}}})$ for which
\[ \gamma \neq 0 \in CH^{2} (A_{\overline{K_{N}}}) \otimes \zit/\ell^{r} \]
is anti-invariant with respect to each $\tilde{\sigma}_{j}$ for $r>>0$. In order to make use of Nori's argument, we need a result that gives a characterization of the minimum field over which $\gamma$ is defined. 
\begin{Def} Given the function field $K =\cit(V)$ of a variety over $\cit$, we say that a finite extension $K \subset L$ is {\em unramified} if there exists an \'etale cover $\tilde{V} \to V$ such that $L = \cit(\tilde{V})$. 
\end{Def}
\begin{Lem}\label{cover} Let $L_{N}$ be an unramified extension of $K_{N}$ and let $Y$ be an \'etale cover $\rho: Y \to W(N)$ such that $L_{N} = \cit(Y)$. Then, there exists \'etale covers $\rho_{j}: X_{j} \to X(N)$ such that \[\rho_{1} \times \rho_{2} \times \rho_{3}: X = X_{1} \times X_{2} \times X_{3} \to W(N)\] factors through $\rho: Y \to W(N)$. 
\begin{proof} We have
\[ \pi_{1}(Y) \leq \pi_{1}(X(N)) \times \pi_{1}(X(N)) \times \pi_{1}(X(N)) = \Gamma(N)^{\times 3}\]
which means that $\pi_{1}(Y)$ is a lattice in the semi-simple Lie group $SL_{2} (\rit)^{\times 3}$. By \cite{WM} Proposition 4.3.3 (a consequence of the Borel Density theorem), there exist lattices $\Gamma_{j} \subset SL_{2} (\rit)$ for $j=1,2,3$ such that $\Gamma = \Gamma_{1} \times \Gamma_{2} \times \Gamma_{3}$ is of finite index in $\pi_{1}(Y)$. It follows that there are finite \'etale covers $\rho_{j}: X_{j} \to X(N)$ for $j=1,2,3$ such that $\pi_{1}(X_{j}) = \Gamma_{j}$. Since
\[ \Gamma \leq \pi_{1}(Y)\leq \pi_{1} (W(N)) \]
there is a finite \'etale cover $X \to Y$ for which the composition \[ X \to Y \to W(N)\] 
is $\rho_{1} \times \rho_{2} \times \rho_{3}$. 
\end{proof}
\end{Lem}
\begin{Cor}\label{cor-cover} In the notation of the above lemma, let $\A_{Y}:= \A(N) \times_{W(N)} Y$. Then, there is an isomorphism of $X$-schemes,
\[ \A_{Y} \times_{Y} X \cong \E_{1} \times \E_{2} \times \E_{3} \]
where $\E_{j} := \E(N) \times_{X(N)} X_{j}$.
\begin{proof} This follows from the fineness of the moduli space $W(N)$ and the fact that the composition $X \to Y \to W(N)$ coincides $\rho_{1} \times \rho_{2} \times \rho_{3}$.
\end{proof}
\end{Cor}
\begin{Prop}\label{crucial} Let $L_{N}$ be an unramified extension of $K_{N}$. Then, for $r >> 0$ \[\gamma \not\in (CH^{2} (A_{\overline{L_{N}}}) \otimes \zit/\ell^{r})^{Gal(\overline{L_{N}}/L_{N})}\]
\begin{proof} Let $K_{N} \subset K_{N}'$ be a finite extension over which \[\gamma \in CH^{2} (A_{\overline{L_{N}}}) \] is defined (i.e., for which $\gamma$ lies in the image of $CH^{2} (A_{K_{N}'}) \to CH^{2} (A_{\overline{L_{N}}})$) and let $L_{N}':= K_{N}'L_{N}$. Our goal will be to show that 
\[\gamma \not\in (CH^{2} (A_{L_{N}'}) \otimes \zit/\ell^{r})^{Gal(L_{N}'/L_{N})}\]
For this, we consider the cycle class:
\begin{equation} P_{*}\gamma \in H^{1}(Gal(\overline{L_{N}}/L_{N}'), PH^{3} (A_{\overline{L_{N}}}, \zit/\ell^{r}(2)))\label{cycle-class} \end{equation}
 We need to show that (\ref{cycle-class}) is not $Gal(L_{N}'/L_{N})$-invariant. To this end,
let $V_{r} := H^{1} (A_{\overline{L_{N}}}, \zit/\ell^{r})$ and as at the beginning of \S \ref{normal} we have
\[W_{r} := V_{r}^{\otimes 3} = PH^{3} (A_{\overline{L_{N}}}, \zit/\ell^{r})  \]
Using inflation-restriction, there is an exact sequence:
\begin{align*} 0 & \to H^{1} (Gal(L_{N}'/L_{N}), W_{r}^{Gal(\overline{L_{N}}/L_{N}')}) \to H^{1} (Gal(\overline{L_{N}}/L_{N}), W_{r})\\ & \to H^{1} (Gal(\overline{L_{N}}/L_{N}'), W_{r})^{Gal(L_{N}'/L_{N})} \to  H^{2} (Gal(L_{N}'/L_{N}), W_{r}^{Gal(\overline{L_{N}}/L_{N}')}) \to \cdots \end{align*}
\begin{claim} $W_{r}^{Gal(\overline{L_{N}}/L_{N}')} = 0$
\begin{proof}[Proof of claim] Since $Gal(\overline{L_{N}}/L_{N}') \leq Gal(\overline{L_{N}}/K_{N})$ has finite index, it suffices to show that $W_{r}$ is irreducible as a $Gal(\overline{L_{N}}/K_{N})$-module. To this end, let $x \in \W(N)(\cit)$ be very general; then, the problem is reduced to showing that $PH^{3} (A_{x}, \zit/\ell^{r})$ is irreducible as a $\pi_{1} (W(N))$-module. This is done in the Appendix.
\end{proof}
\end{claim} 

\noindent Suppose by way of contradiction that (\ref{normal}) is $Gal(L_{N}'/L_{N})$-invariant. Then, by the previous claim, this means that for $r>>0$ 
\[ P_{*}\gamma \in H^{1}(Gal(\overline{L_{N}}/L_{N}), W_{r}) \cong H^{1}_{\text{\'et}} (W(N), PR^{3}g'_{*}\zit/\ell^{r}(2)) \]
\begin{Lem} $\displaystyle P_{*}\gamma \in \mathop{\lim_{\longleftarrow}}_{r>>0} H^{1}_{\text{\'et}} (Y, PR^{3}g'_{*}\zit/\ell^{r}(2))$ is non-torsion.
\begin{proof} It suffices to show that $P_{*}\gamma \neq 0$ after tensoring with $\qit_{\ell}$. Using the comparison isomorphism with singular cohomology, we reduce to showing that the cycle class:
\[ P_{*}\gamma \neq 0 \in H^{1} (Y, PR^{3}g'_{*}\qit(2)) \]
As in \S 3, it suffices to show that $P_{*}\gamma \neq 0 \in H^{4} (\A_{Y}, \qit(2))$, which will follow from the following claim:
\begin{claim} $P_{*}\gamma \neq 0 \in H^{4} (\A_{Y}, \qit(2))$
\begin{proof}[Proof of Claim] By Proposition \ref{key}, this is true for $\A$; one may then pass to some cover $W \to U$ for which $\A_{W} := \A \times_{U} W$ has the full level-$N$ structure so that (by fineness of $W(N)$) there is a (dominant) morphism $W \to W(N)$ with \[\A_{W} \cong \A(N) \times_{W(N)} W\] Then, we have a dominant morphism $W' := W \times_{W(N)} Y \to Y$ for which \[\A_{W} \times_{W} W' \cong \A_{Y} \times_{Y} W'\] By pulling back $\gamma$ via $\A_{W} \times_{W} W' \to \A_{W}$, it follows that \begin{equation} P_{*}\gamma \neq 0 \in H^{4} (\A_{W} \times_{W} W', \qit(2)) \cong H^{4} (\A_{Y} \times_{Y} W', \qit(2))\label{W'} \end{equation} Thus, if $P_{*}\gamma$ is defined over $\A_{Y}$, its pullback to $\A_{W} \times_{W} W'$ is precisely the cycle in (\ref{W'}).
\end{proof}
\end{claim}
\end{proof}
\end{Lem}
\noindent To obtain a contradiction, observe that by pulling back via $\A_{Y} \times_{Y} X \to \A_{Y}$ and using Corollary \ref{cor-cover}, 
\[ P_{*}\gamma \neq 0 \in H^{4} (\E_{1} \times \E_{2} \times \E_{3}, \zit/\ell^{r}(2)) \]
which is anti-invariant with respect to each of the three involutions $\tilde{\sigma}_{j}$ on $\E_{j}$. However, by the K\"unneth theorem, we have
\[ H^{4} (\E_{1} \times \E_{2} \times \E_{3}, \zit/\ell^{r}) \cong \bigoplus_{i_{1}+i_{2}+i_{3} =4} H^{i_{1}} (\E_{1}, \zit/\ell^{r}) \otimes H^{i_{2}} (\E_{2}, \zit/\ell^{r}) \otimes H^{i_{3}} (\E_{3}, \zit/\ell^{r}) \]
It is easy to see that the anti-invariant part of $H^{k} (\E_{j}, \zit/\ell^{r}) = 0$ unless $k=2$. Thus, none of the summands above contains a cycle which is anti-invariant with respect to all three involutions. Hence, the lemma.
\end{proof}
\label{Prop-last}\end{Prop}
\begin{Rem} Since we are working over $SL_{2} (\zit)$, we note that the lattices $\Gamma_{j}$ in the above proof need not be congruence subgroups.
\end{Rem}
\noindent Now, \[L := \mathop{\lim_{\longrightarrow}}_{N \mid n} K_{n}\]
 We consider the absolute Galois group $G_{L}$, and, following Nori's suggestion at the end of \cite{N}, we notice that to every representation:
\[ G_{L} \xrightarrow{\phi} GL_{n} (\zit/\ell^{r}) \]
we may associate a (possibly empty) subset $R(\phi)$ in the following way. Let $H$ be the kernel of $\phi$; then there is some finite Galois extension $L \subset L'$ whose Galois group is the image of $\phi$. This means for each $m \mid N$, there is a compatible system of finite covers $V(n) \to W(n)$ such that $L'$ is the direct limit of $K'_{n} :=\cit(V(N))$; let $R_{n}(\phi)$ be the branch locus of this cover. The branch locus pulls back via $\mathbb{H}^{3} \to W(n)$ (where $\mathbb{H}^{3}$ is the triple product of the upper half plane), and this pullback does not depend on the $n$ chosen. The result of pulling back is $R(\phi)$. 
\begin{Def} We say that $\phi$ is a ramified representation if $R(\phi) \neq \emptyset$.
\end{Def}
\noindent Now, we consider the $\zit/\ell^{r}[G_{L}]$-submodule of $CH^{2} (A_{\overline{L}})\otimes \zit/\ell^{r}$ generated by a $\beta \in CH^{2} (A_{\overline{L}})$; we denote the module by $\M_{\beta, r}$ and the corresponding representation by $\phi_{\beta}$. Then, we have the following lemma:
\begin{Lem} $\phi_{\gamma}$ is a ramified representation.
\begin{proof}[Proof of lemma] This follows directly from Proposition \ref{crucial}. Indeed, if $R(\phi_{\gamma}) = \emptyset$, then $\gamma \in CH^{2} (A_{\overline{L}})\otimes \zit/\ell^{r}$ would be invariant under $Gal(\overline{L}/F)$, where $F$ is the direct limit of unramified extensions $F_{n}$ of $K_{n}$, which would contradict Proposition \ref{Prop-last}. 
\end{proof}
\end{Lem}

\noindent The remainder of the proof will be a sketch, since very similar arguments can be found in \cite{N} and in \cite{T}. There is an action on the field $L$ by a group $G$ which dominates $GL_{2} (\qit)^{\times 3}$, which can be described as follows. Indeed, for \[g \in (GL_{2} (\qit) \cap M_{2} (\zit))^{\times 3},\] there is the usual action on $\mathbb{H}^{3}$ and this action extends to $\mathbb{H}^{3} \times \cit^{3}$. Then, for each $n$, there is some $m \mid n$ and a commutative diagram:
\[ \begin{CD} \A(n) @>{g}>> \A(m)\\
@VVV @VVV\\
W(n)@>{g}>> W(m)  \end{CD}\]
Taking the direct limit over $n$ gives an action of $g$ on $L$ which lifts to an action on the generic fiber, $A_{L}$ (which plays the role of an isogeny). Using the approach of Nori in \cite{N}, one constructs a group $\mathcal{G} \subset Aut(\overline{L}) \times G$ such that for $\tilde{g} = (\tau,g) \in \mathcal{G}$, there is an action on the Chow group \[\tilde{g}^{*}: CH^{2} (A_{\overline{L}}) \to CH^{2} (A_{\overline{L}})\]
which satisfies the following compatibility with the usual action of $G_{L}$ on the Chow group:
\[ \sigma \tilde{g}^{*} = \tilde{g}^{*}\sigma^{\tau^{-1}} : CH^{2} (A_{\overline{L}}) \to CH^{2} (A_{\overline{L}}) \]
and, so long as $\tilde{g}^{*}: CH^{2} (A_{\overline{L}}) \otimes \zit/\ell^{r} \to CH^{2} (A_{\overline{L}})\otimes \zit/\ell^{r}$ is injective, we have
\[ R(\phi_{\tilde{g}^{*}\gamma}) = g^{-1}R(\phi_{\gamma}) \]
(The injectivity can be achieved by ensuring that none of the denominators of the entries of $g \in G$ are divisible by $\ell$; see the Isogeny lemma in \cite{T}). Now, view $SL_{2} (\qit)$ as a subgroup of $G$ via the diagonal imbedding, and we select $g_{1}, g_{2}, \ldots \in SL_{2} (\qit)$ so that $g_{i}$ represent different cosets modulo $SL_{2}(\zit)$ (and so that the denominators of the entries are not divisible by $\ell$). Then, select corresponding lifts $\tilde{g}_{i} = (\tau_{i}, g_{i})\in \mathcal{G}$. We would like to show that there are infinitely many elements in 
\[ \{ \tilde{g}_{1}^{*}\gamma, \tilde{g}_{2}^{*}\gamma, \ldots \} \]
Since $\tilde{g}_{j}^{*}$ were chosen to be injective, it suffices to show the representations $\phi_{\tilde{g}_{i}^{*}\gamma}$ are all distinct, which can be done by showing that the $R_{j}:= R(\phi_{\tilde{g}_{j}^{*}\gamma})$ are all distinct. To this end, we have the following lemma, as in \cite{N} and \cite{T}:
\begin{Lem} The closed subgroup
\[ G_{j} := \{ \beta \in SL_{2}(\qit) \ | \  \beta(R_{j}) = R_{j} \} \]
is a subgroup of $SL_{2} (\zit)$.
\begin{proof}
Because $R_{j}$ is $SL_{2} (\zit)$-stable, the Lie algebra of $G_{j}$ is stable under the adjoint action $SL_{2} (\zit)$ (and, hence, also that of $SL_{2}(\qit)$). Since $\mathfrak{s}\mathfrak{l}_{2}$ is simple, this forces the Lie algebra of $G_{j}$ to be 0. Using the Borel density theorem (\cite{Bo} Theorem 7), it then follows that this group lies inside $SL_{2}(\zit)$. 
\end{proof}
\end{Lem}
\noindent Thus, $gR_{j} = R_{j} \Rightarrow g \in SL_{2} (\zit)$. Since $g_{j}$ were chosen to have distinct cosets modulo $SL_{2} (\zit)$, the result of Theorem \ref{1} now follows.

\section{Appendix}

\begin{Prop} Let $E_{1}$, $E_{2}$ and $E_{3}$ be complex elliptic curves for which the set of $j$-invariants $\{ j(E_{1}), j(E_{2}), j(E_{3}) \}$ has transcendance degree $3$ over $\overline{\qit}$ (so that the product $A = E_{1} \times E_{2} \times E_{3}$ is very general in the moduli of products of elliptic curves). Also, let
\[\iota_{12}: E_{1} \times E_{2} \hookrightarrow A, \ \iota_{23}: E_{2} \times E_{3} \hookrightarrow A, \ \iota_{13}: E_{1} \times E_{3} \hookrightarrow A \]
be the inclusions. Then, the $\ell$ torsion group:
\[ CH^{2} (A) [\ell] = \bigoplus_{j \neq k} \iota_{jk*}CH^{1} (E_{j} \times E_{k})[\ell] \]
for all primes $\ell$.
\begin{proof} Fix a prime $\ell$ and let $N$ be a positive integer for which $2\ell \nmid N$. Then, there exists a fine moduli space of elliptic curves with full level $N$ structure, $X(N)_{\cit}$; let 
\[ f: \E(N)_{\cit} \to X(N)_{\cit}, \  g: \A(N)_{\cit} \times \E(N)_{\cit} \to W(N)_{\cit}  \]
be the corresponding universal elliptic curve and its triple product. It is well-known that $X(N)$ and $\E(N)_{\cit}$ admit models over $\qit(\zeta_{N})$ which have good reduction modulo $\ell$. In particular, the generic fiber of $g$, which we denote by $A$ can be viewed as an Abelian threefold over $\qit(\zeta_{N}) (W(N)) \subset \qit_{\ell}$ with good ordinary reduction. It follows that the main result of \cite{BE} applies and we deduce that
\begin{equation} N^{1}H^{3}_{\text{\'et}} (A_{\overline{\qit_{\ell}}}, \zit/\ell(2)) \neq  H^{3}_{\text{\'et}} (A_{\overline{\qit_{\ell}}}, \zit/\ell(2))\label{BE} \end{equation}
where $N^{*}$ denotes the usual coniveau filtration on \'etale cohomology and, using the comparison isomorphism between \'etale and singular cohomology, we obtain
\[N^{1}H^{3}_{B} (A_{\cit}, \zit/\ell) \neq  H^{3}_{B} (A_{\cit}, \zit/\ell) \]
Now, we observe that
\begin{equation} H^{3}_{B} (A_{\cit}, \zit/\ell) \cong \bigwedge^{3} H^{1}_{B} (A_{\cit}, \zit/\ell)\cong V_{1}^{\oplus 2} \oplus V_{2}^{\oplus 2} \oplus V_{3}^{\oplus 2} \oplus V_{1}\otimes V_{2} \otimes V_{3}\label{decomp} \end{equation} 
where $V_{i} = H^{1} (E_{i,\cit}, \zit/\ell)$ in the notation of the statement of the proposition. We consider the monodromy action of the fundamental group 
\[\pi_{1} (W(N)) = \pi_{1} (X(N)) \times \pi_{1} (X(N))  \times \pi_{1} (X(N)) = \Gamma(N) \times \Gamma(N) \times \Gamma(N)\]
on (\ref{decomp}), for which the $i^{th}$ $\Gamma(N)$ factor of $\pi_{1} (W(N))$ acts on $V_{i}$ via (the reduction modulo $\ell$ of) the standard representation. Since $\ell \nmid N$, the reduction map of $\Gamma(N) \hookrightarrow SL_{2} (\zit) \to SL_{2} (\zit/\ell)$ is surjective. Thus, $\pi_{1} (X(N))$ acts on $V_{i}$ as $SL_{2} (\zit/\ell)$ and, hence, the $V_{1}\otimes V_{2} \otimes V_{3}$ component of (\ref{decomp}) is an irreducible $\pi_{1} (W(N))$-module.  Moreover, the first six summands of $H^{3}_{B} (A_{\cit}, \zit/\ell)$ in (\ref{decomp}) lie in $N^{1}H^{3}_{B} (A_{\cit}, \zit/\ell)$; in fact, first six summands correspond to the image of 
\[ \bigoplus_{jk} H^{1} (E_{j, \cit} \times E_{k, \cit},\zit/\ell) \xrightarrow{\oplus \iota_{jk*}} H^{3}_{B} (A_{\cit}, \zit/\ell) \]
Now, $N^{1}H^{3}_{B} (A_{\cit}, \zit/\ell)$ is invariant under the monodromy action of $\pi_{1}(U)$ for $U \subset W(N)$ open. Since the natural map $\pi_{1} (U) \to \pi_{1} (X)$ is surjective, we may view $H^{3}_{B} (A_{\cit}, \zit/\ell)$ as a $\pi_{1}(U)$-module, and the action of $\pi_{1}(U)$ induces the action of $\pi_{1} (W(N))$ on $H^{3}_{B} (A_{\cit}, \zit/\ell)$. As $\pi_{1}(U)$-modules, we have $N^{1}H^{3}_{B} (A_{\cit}, \zit/\ell) \neq H^{3}_{B} (A_{\cit}, \zit/\ell)$ and $V_{1}\otimes V_{2} \otimes V_{3}$ irreducible. Thus, it follows that 
\[N^{1}H^{3}_{B} (A_{\cit}, \zit/\ell) = V_{1}^{\oplus 2} \oplus V_{2}^{\oplus 2} \oplus V_{3}^{\oplus 2} \]
We deduce from \cite{Su} that Bloch's cycle class map
\[ CH^{2}(A_{\cit})[\ell] \to  H^{3}_{B} (A_{\cit}, \zit/\ell)\]
is an isomorphism whose image is $V_{1}^{\oplus 2} \oplus V_{2}^{\oplus 2} \oplus V_{3}^{\oplus 2}$. However, these summands correspond to the image of 
\[ \bigoplus_{jk} CH^{1} (E_{j, \cit} \times E_{k, \cit})[\ell] \xrightarrow{\oplus \iota_{jk*}} CH^{2}(A_{\cit})[\ell] \]
\end{proof} 
\end{Prop}

\end{document}